\documentclass[12pt]{amsart}
\usepackage{TDefn}

\usepackage[backend=bibtex,  maxbibnames=50]{biblatex}
\newcommand{\triv}{\mathrm{triv}}
\usepackage{cleveref}

\title{Hilbert series of representations of categories of $G$-sets}
\author{Philip Tosteson}
\address{Department of Mathematics, University of North Carolina , Chapel Hill, NC}
\email{\href{mailto:ptoste@unc.edu}{ptoste@unc.edu}}

\begin{document}

\begin{abstract} Let $G$ be a finite group.  A contravariant functor from the category of finite free $G$-sets to vector spaces has an associated Hilbert series, which records the underlying sequence of $G^n$ representations, $n \in \mathbb N$.   We prove that this Hilbert series is rational with denominator given by linear polynomials with coefficients in the field generated by the character table of $G$.  
 \end{abstract}

\maketitle
\section{Introduction}

  Recently, it has been observed that many naturally occurring sequences of group representations arising in algebra and topology admit extra structure: they assemble into a representation of category of a combinatorial nature \cite{CEF, putman2017representation}.   This has motivated the study of the representation theory of combinatorial categories, with a view towards determining the consequences of this additional structure.   In this paper we will consider the particular case of the category of finite $G$-sets.

   Let $G$ be a finite group and let $\FS_G$ be the category of finite free $G$-sets and surjections between them. (Our results will also apply to the category of all maps). Let $k$ be an algebraically closed field of characteristic $p$ (including the case $p = 0$).    An $(\FS_G)\op$ module (or representation) is a functor $M:(\FS_G)\op \to \Vec_k$.   From $M$,  we obtain a sequence of $G^n$ representations  $M_n:= M([n] \times G)$ for $n \in \bbN$,  by evaluating $M$ on the free $G$ set on $[n]$ elements.  Our motivating question is:

\begin{itemize}
 	\item Which sequences of $G^n$ representations can arise in this way?
 \end{itemize}
Actually, in this generality, any sequence is possible.  But when $M$ is \textbf{finitely generated} (see Definition \ref{fingendef}), there are strong restrictions on the possible sequence of representations.  In this paper, we constrain the possible Hilbert series underlying a finitely generated $(\FS_G)\op$ module,  answering a question of Sam--Snowden.   
 
\subsection{Hilbert series and main result}   Following Sam--Snowden \cite[\S6.4]{SamSnowdenG} we define a Hilbert series that determines the class of $M_n$ as a $G^n$ representation in the (rationalized) Grothendieck group.  (Equivalently, determines the Brauer character of $M$).    

 We let $\hat G$ be the set of simple representations of $G$.  We let $\cR(G)$ be the Grothendieck group of $G$ tensored with $\bC$.   Given  $S \in \hat G$ we write $x_S \in \cR(G)$ for the class $S$ so that additively $\cR(G)$ has a basis $\{x_S\}_{S \in \hat H}$.     Then $\cR(G^n) = \cR(G)^{\otimes n} = \bC\{x_S\}_{S \in \hat G}^{\otimes n}$.

 There is a projection map $$\pi: \cR(G^n) = \cR(G)^{\otimes n} \to \Sym^n(\cR(G)),$$ where we take $\Sym^n$ be defined by $\rS_n$ coinvariants. 
Sam--Snowden defined the \textbf{enhanced Hilbert series} of $M$ to be the element $$ H_M = \sum_{n \in \bbN} \pi([M_n])  \in  \prod_{n \in \bbN} \Sym^n(\cR(G)) = \bC[[ x_S ~|~ S \in \hat G]].$$

\begin{ex}\label{freedeg1}
	The free $\FS_G \op$ module on a generator in degree one has $(kG^n)_{n \in \bbN}$ as its underlying sequence of representations.  Its Hilbert series is $$\sum_{n \in \bbN} [kG^n] = \left(\sum_{S \in \hat G} \dim (P_S) x_S \right)^n = \frac{1}{1- \sum_{S \in \hat G} \dim (P_S) x_S},$$ where $P_S$ is the minimal projective cover of $S$.
\end{ex}

\begin{ex}\label{trivial}
	The trivial $\FS_G$ module has $(k\triv_n)_{n \in \bbN}$ as its underlying sequence of representations.  Its Hilbert series is $$\sum_{n \in \bbN} x_{\triv}^n =  \frac{1}{1- x_\triv},$$ where $\triv$ denotes the trivial representation.
\end{ex}

To state our main result, we introduce some notation.  We fix an embedding from the roots of unity of $k$ to $\bC$,  and write $\phi_S$ for the Brauer character of $S$ and $\Phi_S$ for the Brauer character of the minimal projective cover of $S$. 
   Given a $p$-regular conjugacy class $c$ of $G$, let $x_c := \frac{|c|}{|G|} \sum_{S\in \hat G} \overline{ \Phi_S(c) }x_S$.

\begin{thm}\label{mainrationality}
	Let $M$ be a subquotient of an $\FS_G$ module generated in degree $\leq d$.  Then $H_M$ is rational with denominator a product of degree $1$ polynomials of the form $$ 1 - \sum_{c} a_c x_c,$$ where the sum is over $p$-regular conjugacy classes $c$ of $G$,  and  $a_c \in \bbN, a_c \leq |G|d.$ 
 \end{thm}	
 
Both Example \ref{freedeg1} and \ref{trivial} are generated in degree $1$.  These have denominator $1- |G| x_e$ and $1- \sum_{c} x_c$ respectively.  
 
 \begin{remark}
 	Sam--Snowden \cite{SamSnowdenG} proved that $H_M$ was rational with denominator of the form: $1 + \sum_{S} b_S x_S$  for $b_S$ algebraic integers in the field $\bbQ(\zeta_N)$ where $N$ is the exponent of $G$ (or more generally the smallest number such that $G$ is $N$-good).  They asked whether it was possible to improve this class of denominators,  in particular by replacing $\bbQ(\zeta_N)$ by the field generated by the Brauer character table of $\hat G$.   
	Theorem \ref{mainrationality} answers this question in the affirmative, and further restricts the form of the denominators.   \end{remark}

	\subsection{Proof Strategy} We prove Theorem \ref{mainrationality} using the results of our previous paper \cite{TostCat} on $\FSop$ modules.  A finitely generated $(\FS_G)\op$ module restricts to a finitely generated $\FSop$ module, along the functor $X \mapsto G \times X$.  In \cite{TostCat}, for every $\FSop$ module $M$ and $d \in \bbN$ we constructed a sort of Koszul complex $\bK_d(M)$.  Further we proved that if $M$ is finitely generated then, many of the iterated Koszul complexes $$\bK_{d_1} \circ \dots  \circ \bK_{d_r}(M)$$ have vanishing cohomology.   Because the Hilbert series of $\bK_{d}(M)$ is related to the Hilbert series of $M$ in a predictable way,  this vanishing implies that the Hilbert series $H_M$ satisfies nontrivial relations.  In terms of an associated exponential generating function $E_M$,  these relations take the form of a system of linear differential equations.  We diagonalize and solve this system of differential equations.  Because $E_M$ lies in the space of solutions, we obtain the form of $E_M$ and consequently $H_M$.
	
	\subsection{Notation}

\begin{itemize}
	\item $[n]$ denotes the set $\{1, 2, \dots, n\}$.  
	\item $k$ is an algebraically closed field of characteristic $p \geq 0$, $\bC$ is the complex numbers
	\item $G$ is a finite group,  we will also write $G$ for the category with one object with automorphism group $G$
	\item $\hat G$ is the set of isomorphism classes of simple $k$ representations.
	\item $\FS_G$ denotes the category of finite free $G$ sets and surjections  
	\item $\sqcup_n G^n$ denotes the category with objects natural numbers $n \in \bbN$,  and automorphism group $G^n$. 
	\item Given a category $\cC$, we write $\Rep(\cC)$ for the category of functors $\cC \to \Vec_k$, from $\cC$ to the category of $k$ vector spaces. 
	\item $\phi_S$ is the Brauer character of a simple module $S$ and $\Phi_S$ is the Brauer character of its projective cover.
	\item $\cR(G)$ is the Grothendieck group of $G$ tensored with $\bC$, for a $G$ representation $[M]$ denotes its class in $\cR(G)$.  If $S$ is a simple representation, we write $x_S = [S]$.
	\item   $\tr(g,-): \cR(G) \to \bC$ is the Brauer trace, for $g \in G$ a $p$-regular element.
	
\end{itemize}

\section{Variants of Hilbert series and differential equations} There are two other generating functions which record the same data as $H_M$. We define the exponential generating function $E_M \in \bC[[ x_S ~|~ S \in \hat G]]$ by $$E_M := \sum_{n \in \bbN} \frac{\pi([M_n])}{n!},$$ so that $E_M$ is obtained from $H_M$ by applying the transform $$\prod_{S} x_S^{i_S} \mapsto \frac{ \prod_{S} x_S^{i_S}}{(\sum_{S} i_S)!}.$$  Similarly we can define $\tilde H_M$ by applying the transform $\prod_{S} x_S^{i_S} \mapsto \prod_{S} \frac{ x_S^{i_S}}{ i_S!}$.   

For these Hilbert series, we will establish the following version of Theorem \ref{mainrationality}.   Given a $p$-regular conjugacy class $c$ of $G$, let $x_c := \frac{|c|}{|G|} \sum_{S\in \hat G} \overline{ \Phi_S(c) }x_S$.

\begin{thm}\label{exponential}
	Let $M$ be a subquotient of an $(\FS_G)\op$ module generated in degree $\leq d$.  Then 
	\begin{enumerate}
	\item  $E_M$ is a linear combination of products of the form $\prod_{c} (x_c)^{r_c} \exp( a_c x_c)$, for $r_c,a_c \in \bbN$ and $a_c \leq d|G|$.

	\item $\tilde H_M$ is rational, with denominator a product of polynomials of the form $$ P_{c,a}(t) :=\prod_{S \in \hat G} (1-a ~ \frac{|c| ~\overline{\Phi_S(c)}}{|G|} x_S )$$ for $a\in \bbN$, $a \leq d|G|$. 
		\end{enumerate}
 \end{thm}	 
  	
	We note that Theorem \ref{mainrationality} and part (2) of Theorem \ref{exponential} are immediate consequences of part (1), and the fact that the transform $x^n \mapsto n! x^n$ of formal power series acts by $$x^r\exp(jx) \mapsto \left(\frac{d}{dx}\right)^r\frac{1}{1-jx}$$ which is proportional to $\frac{1}{(1-jx)^r}$.

To establish Theorem \ref{mainrationality} we will show that $E_M$ satisfies a system of differential equations.  Given a $p$-regular element $g \in G$, we let $\del_g$ denote the differential operator  $$\del_g:=\sum_{S \in \hat G} \phi_S(g)  \del_S.$$  Note that $\del_g$ only depends on the conjugacy class of $g$.  

\begin{thm}\label{diffequations}
	Let $M$ be an $(\FS_G)\op$ module that is a subquotient of one generated in degree $d$. Then there exists an $r \in \bbN$ such that  $E_M$ satisfies the differential equation ${\del_g \choose d |G| + 1 }^r E_M = 0$ for every $p$-regular $g \in G$.
\end{thm}

Assuming Theorem \ref{diffequations} we now prove part (1) of  Theorem \ref{mainrationality}.  Given a $p$-regular conjugacy class $c$, we have

$$\del_g (x_c) = \frac{|c|}{|G|} \sum_{S \in \hat G} \phi_S(g) \overline{\Phi_S(c)}  = \begin{cases} 1 & \text{ if $g \in c$}  \\ 0  & \text{ otherwise}\end{cases}$$ by the orthogonality relations between simple and projective Brauer characters.   In single variable power series, we have that the solution space of $\prod_{a = 0}^{d |G|} (\del_x - a)^{r}$ is spanned by $x^i \exp(a x)$ for  $a \leq d|G|$ and $i < r$,  $a,i \in \bbN$.  From the orthogonality relation, we obtain that the solution space of differential equations of Theorem \ref{diffequations} is precisely the span of the functions in part (1) of Theorem \ref{mainrationality}.

\section{$\FSop$ modules and Koszul complexes}
In this section we recall the Koszul complexes introduced in \cite{TostCat} and use them to establish Theorem \ref{diffequations}.

Given an $\FSop$ module, $M$, and $d\in \bbN, d \geq 1$ there is an associated Koszul complex $\bK_d(M)$ of $\FSop$ modules \cite[\S 3.2]{TostCat}.  We let $\Sigma^k M$ be the $\FSop$ module defined by $\Sigma^kM(X) := M([k] \sqcup X)$. Then the complex $\bK_d(M)$ takes the form $$\Sigma^d M \leftarrow \Sigma^{d-1} M^{\oplus {d \choose 2}}  \leftarrow \dots \leftarrow  \Sigma^k M^{\oplus s(d,k)}  \leftarrow \dots \leftarrow \Sigma M^{\oplus (d-1)!},$$ where $s(d,k)$ is the unsigned Stirling number of the first kind which is the coefficient of $x^k$ in $\prod_{i = 0}^{d-1} (1 + ix)$.  

We will not fully describe the differentials of $\bK_d(M)$ here; for our purposes, it suffices to know the following.   Given  a surjection  $f: [k] \onto [k-1]$ there is an associated a map of $\FSop$ modules $f^*: \Sigma^{k-1}M \to \Sigma^{k} M$ given by $$(f \sqcup \id_X)^*: M([k-1] \sqcup X) \to M([k] \sqcup X).$$  The differentials of the complex of $\bK_d(M)$ are linear combinations of maps of this form.  

We may iterate the construction of $\bK_d$, by taking the total complex.   The following theorem asserts that for a finitely generated $\FSop$ module, the result of this iteration is often exact.   Before stating it, we recall the definition of finite generation for $\FS_G \op$ modules.  
\begin{defn}\label{fingendef}
	An $\FS_G\op$ module $M$ is finitely generated in degree $\leq d$ if one of the following equivalent conditions holds.
	\begin{enumerate}
		\item There are elements $x_{1} \in M([n_1] \times G),  \dots, x_{r} \in M([n_r] \times G)$ with $n_i \leq d$ such that every $\FS_G \op$ submodule containing $x_1, \dots, x_r$ is equal to $M$.
		\item There is a surjection of $\FS_G \op$ modules $\bigoplus_{i = 1}^r \bP_{\FS_G,n_i} \to M$ where  $n_i \leq d$.
and $\bP_{\FS_G, n_i}$ is the principal projective $\FS_G \op$ module defined by $\bP_{\FS_G, n_i}(X) := k \FS_G(X, [n] \times G)$ 	\end{enumerate}
\end{defn}

\begin{thm}\cite[Theorem 1.2]{TostCat} \label{vanishingthm}
	Let $M$ be an $\FSop$ module which is a subquotient of one that is finitely generated in degree $\leq d$.  Then there exists $r \in \bbN$ such that $\bK_{d+1}^{\circ{r}}(M)$  is exact.  
\end{thm}

\subsection{Restricting $(\FS_G)\op$ modules}  We write $i$ for the embedding of $\FS$ into $\FS_G$ given by $i([n]) = [n] \times G$.  There is a restriction functor $i^*: \Rep(\FS_G \op ) \to \Rep(\FS \op)$.  

There is a functor $G^{-}:  \FS \op \to \Set$ given by $[n] \mapsto G^n$.  If $M$ is an $\FS_G \op$ module, there is a natural action $G^{-} \times i^* M \to i^*M$   given by the action of multiplication $$G^n \times M([n] \times G) \to M([n] \times G).$$
From this construction, we obtain the following.  

\begin{prop}
	The category of $\FS_G \op$ modules is equivalent to the category of $\FS \op$ modules equipped with an action the $\FS\op$ group $G^{-}$ by natural transformations.
\end{prop}
\begin{proof}
	This is a straightforward consequence of the fact that every map of $G$-sets  $[n] \times G \to [m] \times G$ factors uniquely as multiplication by $G^n$ followed by a map of the form $f \times \id$ where $f: [m] \to [n]$.
\end{proof}

Furthermore, we have that finitely generated $\FS_G \op$ modules restrict to finitely generated $\FS\op$ modules.

\begin{prop}\label{fingenrestriction}
	Let $M$ be an $\FS_G \op$ module. If $M$ is a subquotient of a module generated in degree  $\leq d$, then $i^* M$ is of a subquotient of one generated in degree $\leq |G| d$.
\end{prop}
\begin{proof}
	We write $\bP_{\FS_G \op}(n)$ for the principal $\FS_G \op$ module generated in degree $n$.  An $\FS_G \op$ module is The theorem follows from the identity $i^* \bP_{\FS_G \op , n} = \bP_{\FS \op, |G|n}$ and passing to quotients.
\end{proof}

Given an $\FS_G \op$ module, we have that $\bK_d(i^*M)$ carries the structure of a $G \times \FS_G \op$ module, where $G$ acts on $$\Sigma^k i^*M (X) = M([k] \times G \sqcup X \times G)$$ by the constant action on $[k] \times G$ and $G^X $ acts on $\Sigma^k i^*M(X)$ by the action of $G^X$ on $X \times G$.  (The differentials of $\bK_d$ preserve this structure  because the maps $f^*: \Sigma^{k-1} i^* M \to \Sigma^{k} i^* M$ associated above to a surjection $f: [k] \onto [k-1]$ do).  

\subsection{Proof of Theorem \ref{diffequations}}   Forgetting the differentials,  the restriction of $\bK_d(M)$ to $G \times \sqcup_n G^n $ only depends on the restriction of $M$ to $\sqcup_n G^n$.    Since functor $M \mapsto \bK_d(M)$ is exact we obtain an associated operator $$K_d: \prod_{n \in \bbN} \cR(G^n)^{S_n}  \to  \prod_{n \in \bbN}  \cR(G) \otimes \cR(G^n)^{S_n},$$ which takes $[M]$ to the class of $[\bK_d(M)]$. (The class of a finite complex $C_\bdot$ is defined to be $\sum_{i} (-1)^i [C_i] = \sum_{i} (-1)^i H_i(C_\bdot)$).      Further, pairing with the trace operator $\tr(g,- ): \cR(G) \to \bC$ and conjugating by the isomorphism $$\pi/n!: \prod_{n \in \bbN} \cR(G^n)^{S_n} \iso \bC [[x_S | S \in \hat G]]$$ given by projecting and dividing by $n!$  we obtain an operator $$\tr(g,K_d): \bC [[x_S | S \in \hat G]] \to \bC [[x_S | S \in \hat G]].$$

Similiarly there is an operator $\tr(g, \Sigma^k): \bC[[x_S ~|~ S \in \hat G] \to \bC[[x_S ~|~ S \in \hat G] $ for $k \in \bbN$ associated to the functor $\Sigma^k: \Rep(\FS_G \op) \to \Rep(\FS_G \op)$  and we have that $\tr(g, \Sigma^k) = \tr(g, \Sigma)^{\circ k}$

\begin{prop}\label{calc}
	There is an identity of operators $\tr(g,K_d)= \prod_{j = 0}^{d- 1} (\del_g - j)$.
\end{prop}
\begin{proof}
	In homological degree $i$, $\bK_d(M) = \Sigma^{s(d,i)} M$, where $s(d,i)$ is Stirling number defined to be the coefficient of $x^d$ in $\prod_{j = 0}^{d-1}(x + j)$.  We have that $\tr(g, \Sigma^k) = \tr(g, \Sigma)^{\circ k}$, so it suffices to prove that $\tr(g, \Sigma) = \del_g$.  We have that $$\tr(g,\Sigma) x_{S_1} \dots x_{S_m} = \frac{1}{(m-1)!} \pi(\tr(g, \Sigma) \sum_{\sigma \in \bS_m} x_{S_{\sigma(1)}}  \otimes \dots \otimes x_{S_{\sigma(m)}})  $$ $$= \frac{1}{(m-1)!}  \sum_{\sigma \in \bS_m} \phi_{S_{\sigma(1)}}(g)  x_{\sigma(2)} \dots x_{\sigma(m)}. $$ Here $\pi$ denotes the projection $\cR(G)^{\otimes n} \to \Sym^n(\cR(G))$.    This agrees with the action of $\sum_{S \in \hat G} \phi_{S}(g) \del_S$. 
\end{proof}

Given the calculation of Proposition \ref{calc} ,  Theorem \ref{diffequations} is an immediate consequence of the following.  

\begin{prop}
	If $M$ is an $\FS_G$ module which is a subquotient of one generated in degree $d$, exists an $r$ such that  $\rK_{d|G|+1}^{\circ r}(E_M) = 0$.  
\end{prop}
\begin{proof}
	Apply Proposition \ref{fingenrestriction}, and Theorem \ref{vanishingthm}.  Because the homology of $\bK_{d |G| + 1}^{\circ r}(i^* M)$ vanishes,  it follows that the class $\rK_{d|G| + 1}^{\circ r}(E_M) = 0$.  
\end{proof}

\section{Further examples and computations}

In this section,  we include some additional examples of Hilbert series of $\FS_G \op$ modules.

\begin{ex}
	Let $W$ be a projective representation of $G$.  Then there is a $\FS_G \op$ representation $\bQ_W$ with $\bQ_W(X) = \Hom_{G}(W, k\Fin_G(X,G))$.  Then $$E_{\rP_W} = \sum_{c} \frac{|c| \overline{\chi_W(c)}}{|G|} \exp(|G| x_c)$$ and $$H_{P_W} = \sum_{c} \frac{|c| \overline{\chi_W(c)}}{|G|} \frac{1}{1- |G| x_c},$$ where $\chi_W$ denotes the Brauer character of $W$ and the sum is over $p$-regular conjugacy classes of $G$.  
\end{ex}

There is a Day convolution tensor product $\oast$ on the category of $\FS_G \op$ modules defined by $$(M \oast N)(X) = \bigoplus_{X = A \sqcup B} M(A) \otimes N(B).$$ Then $E_{M \oast N} = E_M E_N$.  Combined with the previous example, for any $p$-regular conjugacy class $c$ and any $d \in \bbN$ we may construct an $\FS_G \op$ modules $L$ generated in degree $d$ such that $\exp(d|G| x_c)$ appears in the expansion of $E_L$ with nonzero coefficient.  

\begin{ex}
	Let $M$ be a finitely generated $\FS\op$ module.  Then $M$ pulls back to a finitely generated $(\FS_G)\op$ module along the map $p: \FS_G \to \FS$ defined by $p(X) = X/G$.  The enhanced Hilbert series of $p^*M$ is given by substituting $x_\triv$ into the ordinary Hilbert series of $M$, hence $H_{p^* M}$ has denominator of the form $\prod_{a = 1}^r (1 - a x_{\triv})^{j_a}$.
  \end{ex}

	Finally, we record Proposition \ref{orthogonality} we have that $x_S = \sum_{c} \phi_S(c) x_c$ or more generally that $[M] = \sum_{c} \tr(c,M) x_c$.
	 
	\begin{prop}\label{orthogonality}
		Let $g$ be $p$-regular element of $G$.  Then $$\tr(g, x_c) = \begin{cases} 1 & \text{ if $g \in c$}  \\ 0  & \text{ otherwise}. \\ \end{cases}$$
	\end{prop}
	\begin{proof} 
		This follows from the orthogonality relations between the Brauer characters of simple and projective modules $$\sum_{S \in \hat G} \overline{\Phi_S(g)} \phi_S(h) = \frac{|G|}{|C(g)|} \delta_{g,h}$$
	\end{proof}

\printbibliography

\end{document}